\documentclass[11pt]{article}

\usepackage[T1]{fontenc}
\usepackage[utf8]{inputenc}
\usepackage{lmodern}
\usepackage{microtype}

\usepackage[left=1.05in,right=1.05in,top=1in,bottom=1in]{geometry}
\usepackage{amsmath,amssymb,amsfonts,amsthm,mathtools}
\usepackage{graphicx}
\usepackage{float}
\usepackage{tikz}
\usetikzlibrary{arrows.meta,calc,positioning,decorations.markings}
\usepackage{hyperref}
\usepackage[noabbrev]{cleveref}

\hypersetup{
  colorlinks=true,
  linkcolor=blue!50!black,
  citecolor=blue!60!black,
  urlcolor=blue!60!black
}

\numberwithin{equation}{section}

\usepackage{amsmath,amssymb,amsfonts,amsthm}
\usepackage{array}
\usepackage{graphicx}
\usepackage{mathrsfs}
\usepackage{enumitem}
\usepackage{hyperref}
\usepackage[noabbrev]{cleveref}
\usepackage[nottoc,notlot,notlof]{tocbibind}

\usepackage{booktabs}
\usepackage{tabularx}
\usepackage{array}

\usepackage{tikz-cd}

\newtheorem{theorem}{Theorem}[section]
\newtheorem{proposition}[theorem]{Proposition}

\theoremstyle{definition}
\newtheorem{definition}[theorem]{Definition}
\theoremstyle{remark}

\title{Gauss Maps in Hyperbolic Surface Theory:\\
A Unified Perspective}
\author{Magdalena Toda and Erhan G\"uler}
\date{}

\begin{document}
\maketitle

\begin{abstract}
Immersed surfaces in hyperbolic three-space carry several natural
Gauss-type maps with distinct geometric roles. The hyperbolic Gauss maps
record the ideal endpoints of oriented normal geodesics; the Legendre Gauss
lift retains the position-normal data and its contact structure; adjusted
Gauss maps arise from gauge normalization and Iwasawa splitting in
Weierstrass--Kenmotsu representations; and the conformal Gauss map encodes
the mean-curvature sphere congruence in Möbius geometry. We present these
constructions in a common framework, emphasizing their target spaces,
analytic properties, and mutual relations. Particular attention is given to
the generalized DPW method, the necessity of flatness in adjusted
rank-one data, and the harmonic-map characterization of Willmore surfaces.
The resulting viewpoint distinguishes the asymptotic, contact, integrable,
and conformal information carried by an immersed surface in
\(\mathbb{H}^{3}(-1)\).
\\[2pt]
\end{abstract}

\tableofcontents

\section{Introduction}

Minimal and constant mean curvature surfaces occupy a central place in
differential geometry, geometric analysis, and integrable systems. A guiding
principle of the subject is that the geometry of an immersion can often be
encoded by holomorphic or harmonic data. In Euclidean three-space, the
classical Weierstrass representation reconstructs a minimal surface from its
meromorphic Gauss map and a holomorphic differential \cite{Osserman}. More
generally, Ruh and Vilms proved that the Gauss map of an immersed submanifold
is harmonic if and only if its mean curvature vector is parallel in the normal
bundle \cite{RuhVilms}. For oriented surfaces in $\mathbb{E}^{3}$, this
condition is equivalent to constancy of the mean curvature.

Important global developments further illustrate the richness of this theory.
Lawson constructed compact embedded minimal surfaces of arbitrary genus in
$\mathbb{S}^{3}$ \cite{Lawson}. Pinkall and Sterling showed that constant mean
curvature tori are of finite type and developed their classification through
algebro-geometric spectral data \cite{PinkallSterling}. The broader
differential-geometric background of these developments, including
minimal-submanifold theory and global geometric analysis, is represented in
the collection edited by Yau \cite{Yau}.

This correspondence provides a natural entry point to integrable-systems
methods. Harmonic maps into symmetric spaces admit associated families of
flat connections, and the DPW construction reconstructs them from
holomorphic potentials by means of loop-group factorization \cite{DPW}.
Pressley and Segal established the analytic and geometric foundations of
loop-group theory that underlie many of these constructions
\cite{PressleySegal}. Bobenko formulated classical and integrable surface
geometries in terms of $2\times2$ matrix representations \cite{Bobenko},
while Bobenko and Eitner developed the relationship between special surface
geometries and Painlev\'e equations \cite{BobenkoEitner}. Burstall and Pedit
studied harmonic maps through Adler--Kostant--Symes theory
\cite{BurstallPedit}, whereas Burstall and Guest obtained Weierstrass-type
formulas and classification results for harmonic two-spheres in compact
symmetric spaces \cite{BurstallGuest}. In hyperbolic three-space, however,
the classical roles played by the Euclidean Gauss map separate into several
geometrically distinct constructions. No single Gauss map simultaneously
captures the asymptotic, contact, integrable, and conformal structures
associated with an immersed surface.

The first natural construction is the pair of hyperbolic Gauss maps
\begin{equation*}
G^{\pm}\colon M\longrightarrow\partial_{\infty}\mathbb{H}^{3},
\end{equation*}
which assign to each point the two ideal endpoints of the oriented normal
geodesic. These maps encode the asymptotic geometry of the normal congruence.
Their systematic study is closely associated with the work of Epstein
\cite{Epstein}, and the positive hyperbolic Gauss map plays a central role in
Bryant's representation of CMC-one surfaces in hyperbolic space
\cite{Bryant}. Their global significance for complete CMC-one surfaces was
further developed by Umehara and Yamada \cite{UmeharaYamada}.

A second construction retains the complete position--normal pair
\begin{equation*}
\mathcal{F}=(f,N)\colon M\longrightarrow U\mathbb{H}^{3}.
\end{equation*}
This Legendre Gauss lift satisfies the canonical contact condition and
provides framed data naturally adapted to the $4$-symmetric-space formulation
of hyperbolic CMC surfaces. Dorfmeister, Inoguchi, and Kobayashi used this
framework to derive a generalized Weierstrass representation, with particular
emphasis on the range $0\leq H<1$ \cite{DIK}.

Adjusted Gauss maps arise after an appropriate gauge normalization together
with an Iwasawa factorization. Unlike the preceding constructions, they are
not canonical maps determined solely by the immersion. Rather, they are
representation-dependent objects designed to simplify reconstruction
formulas. In the Weierstrass--Kenmotsu representation of Toda, G\"uler, and
Atampalage, CMC immersions with $0\leq H<1$ are locally described by flat
rank-one data
\begin{equation*}
\Omega=\eta-\lambda\eta^{*},
\qquad
d\Omega+\Omega\wedge\Omega=0,
\end{equation*}
where
\begin{equation*}
H=\frac{1-|\lambda|^{2}}{1+|\lambda|^{2}}.
\end{equation*}
The flatness condition is an essential differential compatibility requirement
and does not follow from the algebraic condition $\det\eta=0$
\cite{TodaGulerAtampalage}.

A fourth construction belongs to Möbius geometry. Within the broader
variational background, Bryant and Griffiths developed a systematic reduction
method for constrained variational problems and applied it to the elastic-curve
functional
\begin{equation*}
\int \frac{1}{2}\kappa^{2}\,ds
\end{equation*}
\cite{BryantGriffiths}. The conformal Gauss map assigns to an immersion its
mean-curvature sphere congruence. Its harmonicity characterizes Willmore
surfaces and leads naturally to a DPW formulation in a conformal symmetric
space. This viewpoint has been developed through the work of H\'elein
\cite{Helein}, Xia and Shen \cite{XiaShen}, Dorfmeister and Wang
\cite{DorfmeisterWang}, and Wang \cite{Wang}.

The purpose of this article is to clarify the geometric role of each of these
constructions and to explain the relationships among them. The hyperbolic
Gauss maps describe geometry at infinity; the Legendre lift retains contact
position--normal data; adjusted Gauss maps organize flat gauge-normalized data
for explicit reconstruction; and the conformal Gauss map records the
Möbius-invariant mean-curvature sphere congruence. These maps should therefore
not be viewed as competing notions of a Gauss map, but rather as complementary
geometric objects that emphasize different aspects of the same immersion.

We begin by recalling the classical Euclidean picture before examining the
hyperbolic, Legendre, adjusted, and conformal Gauss maps in turn. The final
sections compare their target spaces, defining structures, analytic
properties, and geometric significance from a unified perspective.

%%%%%%%%%%%%%%%%%%%%%%%%%%%%%%%%%%%%%%%%%%%%%%%%%%%%%%%%%%%%%%%
%%%%%%%%%%%%%%%%%%%%%%%%%%%%%%%%%%%%%%%%%%%%%%%%%%%%%%%%%%%%%%%

\section{The Classical Euclidean Picture}
Let $f:M\to\mathbb{R}^{3}$ be an oriented immersed surface with unit normal
vector field $N$. The shape operator (or Weingarten map) is defined by
\[
S=-\,dN,
\]
or equivalently,
\[
dN=-df\circ S.
\]
Its eigenvalues are the principal curvatures of the surface. Consequently,
\[
K=\det S,\qquad
H=\frac12\operatorname{tr}S,\qquad
|dN|^{2}=|S|^{2}=4H^{2}-2K.
\]
Thus, the differential of the Gauss map completely encodes the extrinsic
curvature of the immersion.

Equip \(M\) with the conformal structure induced by \(f\). For a conformal minimal immersion, stereographic projection identifies the Gauss map with a meromorphic function
\[
g\colon M\longrightarrow\mathbb{C}\cup\{\infty\}.
\]
Using stereographic projection from the north pole,
\[
g=\frac{N_{1}+iN_{2}}{1-N_{3}}.
\]

Conversely, let \(g\) be meromorphic and let \(\eta\) be a holomorphic \(1\)-form such that
\[
\Phi=
\left(
\frac12(1-g^{2}),
\frac{i}{2}(1+g^{2}),
g
\right)\eta
\]
is holomorphic and has no common zeros. Since
\[
\langle\Phi,\Phi\rangle=0,
\]
the map
\[
f=\operatorname{Re}\int\Phi
\]
defines a conformal minimal immersion on a simply connected domain. On a non-simply connected surface, the period condition
\[
\operatorname{Re}\oint_{\gamma}\Phi=0,
\qquad
\gamma\in H_{1}(M,\mathbb{Z}),
\]
must also be satisfied. The induced metric is
\[
ds^{2}
=
\frac14(1+|g|^{2})^{2}|\eta|^{2}.
\]
Hence, the Weierstrass data \((g,\eta)\) determine the immersion uniquely up to translations of \(\mathbb{E}^{3}\) \cite{Osserman}.

For a general conformal immersion, the Hopf differential is
\[
Q\,dz^{2}
=
\langle f_{zz},N\rangle\,dz^{2}.
\]
The Codazzi equations imply that \(Q\,dz^{2}\) is holomorphic whenever the mean curvature is constant. In the minimal case, and with the above normalization,
\[
Q=-\eta\,dg,
\]
up to the sign convention adopted for the second fundamental form.

A fundamental analytic characterization of the Euclidean Gauss map is provided by the Ruh--Vilms theorem. In its general form, the Gauss map of an immersed Euclidean submanifold is harmonic if and only if its mean curvature vector is parallel in the normal bundle \cite{RuhVilms}. For oriented surfaces in \(\mathbb{E}^{3}\), this becomes the following classical result.

\begin{theorem}[Ruh--Vilms]
Let
\[
f\colon M\longrightarrow\mathbb{E}^{3}
\]
be an oriented immersed surface with Gauss map
\[
N\colon M\longrightarrow\mathbb{S}^{2}.
\]
Then \(N\) is harmonic if and only if the mean curvature \(H\) is constant.
\end{theorem}

Accordingly,
\[
\Delta N+|dN|^{2}N=0,
\]
where \(\Delta\) denotes the Laplace--Beltrami operator. In particular, the Gauss map of a minimal surface is simultaneously harmonic as a map into \(\mathbb{S}^{2}\) and meromorphic after the identification
\[
\mathbb{S}^{2}\cong\mathbb{C}\cup\{\infty\}.
\]
Thus, harmonicity and holomorphicity provide two complementary analytic descriptions of the same geometric object.
Since
\[
\mathbb{S}^{2}\cong SO(3)/SO(2)
\]
is a Riemannian symmetric space, the harmonic map equation admits a remarkable reformulation. Instead of studying the nonlinear second-order differential equation satisfied by the Gauss map directly, one associates to it a family of linear first-order systems depending on an auxiliary complex parameter, called the \emph{spectral parameter}. The harmonicity of the Gauss map is then equivalent to the compatibility, or zero-curvature, condition for this family of linear systems. This reformulation makes available powerful methods from integrable systems, Lie theory, and loop groups. In particular, the DPW method reconstructs harmonic maps from holomorphic potential data by integrating these linear systems and applying loop-group factorizations \cite{DPW}.

The Euclidean theory therefore provides the prototype for the developments that follow. The differential of the Gauss map encodes curvature, meromorphicity yields the classical Weierstrass representation for minimal surfaces, and harmonicity gives rise to an integrable-systems formulation. In hyperbolic geometry, however, these different geometric roles are no longer carried by a single object but instead separate naturally into several distinct Gauss-type constructions.
%%%%%%%%%%%%%%%%%%%%%%%%%%%%%%%%%%%%%%%%%%%%%%%%%%%%%%%%%%%%%%%
%%%%%%%%%%%%%%%%%%%%%%%%%%%%%%%%%%%%%%%%%%%%%%%%%%%%%%%%%%%%%%%
%%%%%%%%%%%%%%%%%%%%%%%%%%%%%%%%%%%%%%%%%%%%%%%%%%%%%%%%%%%%%%%
\section{Hyperbolic Gauss Maps}

The classical Gauss map of an oriented surface in Euclidean
three-space assigns to each point its unit normal vector, regarded
as a point of the fixed sphere \(\mathbb S^{2}\). This construction
uses the translational structure of Euclidean space, which identifies
tangent spaces at distinct points. Hyperbolic space has no analogous
translation structure. Instead, the asymptotic geometry of its
geodesics leads naturally to Gauss maps with values on the ideal
boundary. This point of view was developed systematically by Epstein
in his study of the hyperbolic Gauss map and quasiconformal geometry
\cite{Epstein}.

Throughout this section, we use the hyperboloid model
\[
\mathbb H^{3}(-1)
=
\left\{
x\in\mathbb R^{1,3}:
\langle x,x\rangle_{1,3}=-1,\quad x_{0}>0
\right\},
\]
where
\[
\langle x,y\rangle_{1,3}
=
-x_{0}y_{0}
+x_{1}y_{1}
+x_{2}y_{2}
+x_{3}y_{3}
\]
is the Lorentzian inner product on \(\mathbb R^{1,3}\).

The ideal boundary of hyperbolic space is naturally identified with
the projectivized positive light cone
\[
\partial_{\infty}\mathbb H^{3}
\cong
\mathbb P(\mathcal L^{+}),
\]
where
\[
\mathcal L^{+}
=
\left\{
\xi\in\mathbb R^{1,3}\setminus\{0\}:
\langle\xi,\xi\rangle_{1,3}=0,\quad \xi_{0}>0
\right\}.
\]
Thus, a point of \(\partial_{\infty}\mathbb H^{3}\) is represented by
a positive lightlike vector, with two such vectors representing the
same point when they differ by multiplication by a positive scalar.
The ideal boundary carries a natural conformal structure and is
conformally equivalent to the Riemann sphere.

Let
\[
f:M\longrightarrow\mathbb H^{3}(-1)
\]
be an oriented immersion, and let \(N\) be its unit normal vector
field. Since \(\mathbb H^{3}(-1)\) is the level set
\[
\langle x,x\rangle_{1,3}=-1,
\]
its tangent space at \(f(p)\) is the Lorentzian orthogonal complement
of \(f(p)\):
\[
T_{f(p)}\mathbb H^{3}(-1)
=
\left\{
v\in\mathbb R^{1,3}:
\langle f(p),v\rangle_{1,3}=0
\right\}.
\]
Consequently,
\[
\langle f,f\rangle_{1,3}=-1,\qquad
\langle N,N\rangle_{1,3}=1,\qquad
\langle f,N\rangle_{1,3}=0.
\]
It follows that
\[
\langle f\pm N,f\pm N\rangle_{1,3}=0,
\]
so \(f+N\) and \(f-N\) determine points of the projectivized light
cone.

\begin{definition}
The \emph{positive} and \emph{negative hyperbolic Gauss maps} of the
oriented immersion \(f\) are the maps
\[
G^{\pm}:M\longrightarrow\partial_{\infty}\mathbb H^{3},
\qquad
G^{\pm}(p)=[f(p)\pm N(p)],
\]
where \([\xi]\) denotes the projective class of a nonzero positive
lightlike vector.
\end{definition}

Geometrically, \(G^{+}(p)\) and \(G^{-}(p)\) are the two endpoints at
infinity of the oriented normal geodesic through \(f(p)\). Indeed, that
geodesic is
\[
\gamma_{p}(t)
=
\cosh(t)\,f(p)+\sinh(t)\,N(p).
\]
Since projective classes are unchanged by multiplication by a nonzero
scalar,
\[
[\gamma_{p}(t)]
=
[f(p)+\tanh(t)N(p)].
\]
Therefore,
\[
[\gamma_{p}(t)]
\longrightarrow
[f(p)+N(p)]
\qquad (t\to+\infty),
\]
whereas
\[
[\gamma_{p}(t)]
\longrightarrow
[f(p)-N(p)]
\qquad (t\to-\infty).
\]

Thus, unlike the Euclidean Gauss map, which records a unit normal
direction in a fixed sphere, the hyperbolic Gauss maps record the
asymptotic endpoints of the normal geodesic congruence. Epstein
introduced this construction and developed its relation to the
geometry of surfaces and quasiconformal mappings
\cite{Epstein}. Bryant subsequently showed that the positive
hyperbolic Gauss map is a fundamental holomorphic invariant of
constant mean curvature one surfaces in hyperbolic space
\cite{Bryant}.

Changing the orientation replaces \(N\) by \(-N\), and hence
interchanges \(G^{+}\) and \(G^{-}\). Thus neither map is individually
independent of orientation, although the unordered pair
\[
\{G^{+},G^{-}\}
\]
is intrinsic to the immersed surface.

\subsection{Differential and conformal properties}

Let
\[
S:TM\longrightarrow TM
\]
denote the shape operator, with the convention
\[
dN=-df\circ S.
\]
For the light-cone lifts
\[
Y^{\pm}=f\pm N,
\]
the Weingarten equation gives
\[
dY^{\pm}
=
df\circ(I\mp S).
\]
This elementary formula contains the principal local information
carried by the two hyperbolic Gauss maps.

In particular, \(G^{+}\) is singular at \(p\) precisely when
\(I-S\) is singular at \(p\), whereas \(G^{-}\) is singular precisely
when \(I+S\) is singular. Hence
\[
G^{+}\text{ is singular at }p
\quad\Longleftrightarrow\quad
\det(I-S)_{p}=0,
\]
and
\[
G^{-}\text{ is singular at }p
\quad\Longleftrightarrow\quad
\det(I+S)_{p}=0.
\]
If \(k_{1}\) and \(k_{2}\) are the principal curvatures, these
conditions become
\[
k_{i}=1
\quad\text{for at least one }i
\]
for \(G^{+}\), and
\[
k_{i}=-1
\quad\text{for at least one }i
\]
for \(G^{-}\). These singularity conditions are manifestations of the
contact of the surface with horospheres; see Epstein
\cite{Epstein}.

The conformal structures induced by the light-cone lifts may be
expressed in terms of the three fundamental forms. For tangent vectors
\(X,Y\in TM\),
\[
\langle dY^{+}(X),dY^{+}(Y)\rangle_{1,3}
=
I(X,Y)-2\,\mathrm{II}(X,Y)+\mathrm{III}(X,Y),
\]
and
\[
\langle dY^{-}(X),dY^{-}(Y)\rangle_{1,3}
=
I(X,Y)+2\,\mathrm{II}(X,Y)+\mathrm{III}(X,Y).
\]
Accordingly, we write
\[
I_{\infty}^{+}
=
I-2\,\mathrm{II}+\mathrm{III},
\qquad
I_{\infty}^{-}
=
I+2\,\mathrm{II}+\mathrm{III}.
\]

Here and below, a map from the Riemann surface determined by the
induced metric \(I\) to the conformal sphere
\(\partial_{\infty}\mathbb H^{3}\) is called
\emph{weakly conformal} if its pullback metric has the form
\[
G^{*}g_{\partial_{\infty}\mathbb H^{3}}
=
\lambda I
\]
for some nonnegative function \(\lambda\). Thus the map is conformal
where \(\lambda>0\), while points at which \(\lambda=0\) are allowed
as branch or degenerate points.

Let \(e_{1},e_{2}\) be an orthonormal frame of principal directions.
Then
\[
I_{\infty}^{+}(e_i,e_i)
=
(1-k_i)^2,
\qquad
I_{\infty}^{+}(e_1,e_2)=0,
\]
and
\[
I_{\infty}^{-}(e_i,e_i)
=
(1+k_i)^2,
\qquad
I_{\infty}^{-}(e_1,e_2)=0.
\]
Therefore, \(G^{+}\) is weakly conformal precisely when
\[
(1-k_{1})^{2}=(1-k_{2})^{2}.
\]
This equality is equivalent to
\[
(k_{1}-k_{2})(k_{1}+k_{2}-2)=0.
\]
Hence either the point is umbilic,
\[
k_{1}=k_{2},
\]
or
\[
k_{1}+k_{2}=2.
\]
With the convention
\[
H=\frac{k_{1}+k_{2}}{2},
\]
the latter condition is \(H=1\).

Similarly, \(G^{-}\) is weakly conformal precisely when
\[
(1+k_{1})^{2}=(1+k_{2})^{2},
\]
or equivalently
\[
(k_{1}-k_{2})(k_{1}+k_{2}+2)=0.
\]
Thus, away from umbilic points, weak conformality of \(G^{-}\) is
equivalent to
\[
H=-1.
\]

This computation gives the following local characterization, whose
positive-map formulation is central to Bryant's theory of CMC-one
surfaces \cite{Bryant}; the corresponding statement for \(G^{-}\)
is obtained by reversing the orientation.

\begin{proposition}[Epstein--Bryant characterization]
Let
\[
f:M\longrightarrow\mathbb H^{3}(-1)
\]
be an oriented immersed surface.

\begin{enumerate}[label=\rm(\roman*)]

\item
The positive hyperbolic Gauss map \(G^{+}\) is weakly conformal at
every umbilic point. At every nonumbilic point, it is weakly conformal
if and only if
\[
H=1.
\]

\item
The negative hyperbolic Gauss map \(G^{-}\) is weakly conformal at
every umbilic point. At every nonumbilic point, it is weakly conformal
if and only if
\[
H=-1.
\]

\item
Consequently, if the surface is not totally umbilic, then \(G^{+}\)
is weakly conformal on \(M\) if and only if the immersion has constant
mean curvature \(H=1\). Likewise, \(G^{-}\) is weakly conformal on
\(M\) if and only if the immersion has constant mean curvature
\(H=-1\).

\end{enumerate}
\end{proposition}

In particular, Bryant's result shows that CMC-one surfaces occupy in
hyperbolic geometry a position analogous to that of minimal surfaces
in Euclidean space: the positive hyperbolic Gauss map of a
non-totally-umbilic CMC-one immersion is conformal, or equivalently
meromorphic, with respect to the complex structure induced by the
first fundamental form \cite{Bryant}. For the oppositely oriented
surface, the corresponding distinguished map is \(G^{-}\).

\subsection{Bryant's representation}

The hyperbolic Gauss map assumes a distinguished role in Bryant's
representation of constant mean curvature one surfaces in
hyperbolic three-space \cite{Bryant}. This representation is the
hyperbolic analogue of the classical Weierstrass representation for
minimal surfaces in Euclidean space and expresses a CMC-one immersion
in terms of holomorphic data.

Let \(M\) be simply connected, and let
\[
F:M\longrightarrow SL(2,\mathbb C)
\]
be a holomorphic null immersion. In the Hermitian matrix model of
\(\mathbb H^{3}(-1)\), Bryant showed that
\[
f=FF^{*},
\qquad
F^{*}=\overline{F}^{\,t},
\]
defines a conformal immersion of constant mean curvature one
\cite{Bryant}.

Locally, the Maurer--Cartan form of \(F\) may be written as
\[
F^{-1}dF
=
\begin{pmatrix}
g & -g^{2}\\
1 & -g
\end{pmatrix}
\omega,
\]
where \(g\) is meromorphic and \(\omega\) is a holomorphic
\(1\)-form. The resulting induced metric is
\[
ds^{2}
=
(1+|g|^{2})^{2}|\omega|^{2},
\]
and the Hopf differential is
\[
Q=\omega\,dg,
\]
up to the sign convention used for the second fundamental form
\cite{Bryant}.

Although \(g\) resembles the stereographic Gauss map appearing in the
Euclidean Weierstrass representation, it is not the geometric
hyperbolic Gauss map defined by the endpoints of the normal
geodesics. It is commonly called the \emph{secondary Gauss map}. It
arises from the holomorphic frame and describes the complex analytic
data of Bryant's representation.

By contrast, the geometric hyperbolic Gauss map is determined
directly by the immersion and its oriented normal congruence. If
\[
F=
\begin{pmatrix}
F_{11}&F_{12}\\
F_{21}&F_{22}
\end{pmatrix},
\]
then, wherever the indicated quotients are defined, Bryant's
representation gives
\[
G
=
\frac{dF_{11}}{dF_{21}}
=
\frac{dF_{12}}{dF_{22}}
\]
\cite{Bryant}. Thus a CMC-one immersion naturally gives rise to two distinct
meromorphic maps:
\begin{itemize}
\item the geometric hyperbolic Gauss map \(G\), determined by the
      endpoints at infinity of the oriented normal geodesics; and
\item the secondary Gauss map \(g\), determined by the holomorphic
      frame in Bryant's representation.
\end{itemize}
Although both maps are meromorphic, they encode different geometric
information and should not be identified.

Their relationship is expressed through the Schwarzian derivative.
For a locally univalent meromorphic function \(h\), set
\[
S(h)
=
\left[
\left(\frac{h''}{h'}\right)'
-
\frac12
\left(\frac{h''}{h'}\right)^{2}
\right]dz^{2}.
\]
Then
\[
S(g)-S(G)=2Q,
\]
up to the sign convention adopted for the Hopf differential. Thus the
Hopf differential measures the difference between the complex
projective structures determined by the secondary and geometric
Gauss maps.

The global theory of complete CMC-one surfaces was developed further
by Umehara and Yamada \cite{UmeharaYamada}. In particular, they studied
complete CMC-one surfaces of finite total curvature, their
conformal compactification, and the behavior of the Hopf differential
and the hyperbolic Gauss map at the ends
\cite{UmeharaYamada}.

In this theory, the Hopf differential extends meromorphically across
the punctures arising in the conformal compactification. The
hyperbolic Gauss map, however, need not extend meromorphically across
every end. This distinction gives rise to the notions of
\emph{regular} and \emph{irregular} ends, which play an important role
in the global study of CMC-one surfaces
\cite{UmeharaYamada}.

The hyperbolic Gauss maps record the asymptotic behavior of the
oriented normal geodesics, but they do not retain the full
position--normal data of the immersion or the contact-geometric
structure underlying the integrable-systems formulation. Those
additional structures are incorporated by the Legendre Gauss lift,
which is discussed in the next section.
%%%%%%%%%%%%%%%%%%%%%%%%%%%%%%%%%%%%%%%%%%%%%%%%%%%%%%%%%%%%%%%
%%%%%%%%%%%%%%%%%%%%%%%%%%%%%%%%%%%%%%%%%%%%%%%%%%%%%%%%%%%%%%%

%%%%%%%%%%%%%%%%%%%%%%%%%%%%%%%%%%%%%%%%%%%%%%%%%%%%%%%%%%%%%%%
\section{The Legendre Gauss Lift and the Generalized DPW Method}

The hyperbolic Gauss maps introduced in the previous section record the
two ideal endpoints of the oriented normal geodesic through each point of
an immersed surface. Although these maps capture important asymptotic
information, they do not determine the complete contact element, namely
the point of the surface together with its oriented tangent plane. For
many geometric questions it is therefore natural to retain both the
position and the unit normal simultaneously.

Let
\[
f:M\longrightarrow\mathbb H^{3}(-1)
\]
be an oriented immersion with unit normal vector field \(N\). The map
\[
L=(f,N):M\longrightarrow U\mathbb H^{3}(-1)
\]
is called the \emph{Legendre Gauss lift}. Since the unit tangent bundle
\(U\mathbb H^{3}(-1)\) carries a natural contact structure, the image of
\(L\) is Legendrian. Thus the Legendre lift preserves the complete
position--normal data of the immersion rather than only the asymptotic
information encoded by the hyperbolic Gauss maps.

From the geometric point of view, the Legendre lift is therefore a
natural refinement of the classical Gauss map. It associates to each
point of the surface its oriented contact element and provides a
convenient framework for studying transformations, congruences of
normal geodesics, and deformations of immersed surfaces.

A further development was made by Dorfmeister, Inoguchi, and Kobayashi,
who showed that the unit tangent bundle of hyperbolic three-space admits
a homogeneous-space realization as a \(4\)-symmetric space and that the
Legendre lift of a constant mean curvature surface can be interpreted
within the corresponding harmonic-map theory \cite{DIK}. This
interpretation connects classical surface theory with harmonic maps,
Lie groups, and integrable systems.

One consequence of this viewpoint is that the harmonicity of the
Legendre lift may be encoded by a one-parameter family of flat
connections depending on an auxiliary complex parameter
\(\lambda\), called the \emph{spectral parameter}. Equivalently, one
constructs an extended frame whose Maurer--Cartan form satisfies a
zero-curvature condition for every value of \(\lambda\). In this way,
the nonlinear geometric problem is transformed into an equivalent
family of compatible linear systems.

As in the classical DPW construction for harmonic maps into symmetric
spaces \cite{DPW}, the extended frame is generated from holomorphic
data, commonly called a \emph{holomorphic potential}. Integration of
this potential, followed by suitable loop-group decompositions such as
the Iwasawa and Birkhoff factorizations, reconstructs the harmonic
Legendre lift and hence the corresponding immersed surface. The
analytic details of this generalized DPW method are developed
systematically by Dorfmeister, Inoguchi, and Kobayashi
\cite{DIK}, and will not be repeated here.

An important feature of this approach is its unifying character.
Various Weierstrass-type representation formulas for constant mean
curvature surfaces in hyperbolic space, originally developed in
different geometric settings, may be interpreted within a common
harmonic-map framework. Bryant's representation for CMC-one surfaces
\cite{Bryant}, the generalized representations for surfaces satisfying
\(0\le H<1\), and the corresponding constructions for \(H>1\) can all
be viewed from this broader perspective, although their geometric
motivations and analytic realizations remain distinct.

From the perspective of the present survey, however, the principal role
of the Legendre lift is geometric rather than analytic. It enriches the
hyperbolic Gauss map by recording the complete contact element of the
surface and provides a natural bridge between classical differential
geometry, contact geometry, harmonic maps, and modern representation
theory.

The Legendre lift is one possible extension of the hyperbolic Gauss
map. In the next section we consider a different approach, based on
adjusted Gauss maps, which is motivated more directly by explicit
representation formulas.%%%%%%%%%%%%%%%%%%%%%%%%%%%%%%%%%%%%%%%%%%%%%%%%%%%%%%%%%%%%%%%%%%
\section{Adjusted Gauss Maps and Generalized Weierstrass Representations}

The Legendre Gauss lift provides the natural contact-geometric framework
for the generalized DPW method. For explicit representation formulas,
however, it is advantageous to introduce an additional normalization of
the moving frame. This leads to the \emph{adjusted Gauss map}, a
representation-adapted Gauss map arising naturally from a balanced
matrix formulation of constant mean curvature immersions. Although it is
not a canonical invariant of the immersion, it provides the geometric
normalization that makes generalized Weierstrass--Kenmotsu
representations particularly transparent.

Unlike the hyperbolic Gauss maps and the Legendre lift, whose
definitions depend only on the immersion and its orientation, the
adjusted Gauss map depends on a balanced choice of gauge. Its purpose is
not to introduce additional geometric information, but rather to
reorganize the existing position--normal data into a form well suited
for explicit reconstruction, Iwasawa factorization, and the associated
loop-group formalism.

Let
\[
S:M\longrightarrow SL(2,\mathbb C)
\]
be an adapted moving frame associated with a constant mean curvature
immersion. Following Toda, G\"uler, and Atampalage
\cite{TodaGulerAtampalage}, one performs a balanced gauge
transformation together with an Iwasawa decomposition
\[
S=F_s\Phi,
\]
where \(F_s\) belongs to the noncompact factor and
\(\Phi\in SU(2)\). The unitary factor determines the adjusted Gauss map
\[
N_{\mathrm{ad}}
=
\Phi\sigma_3\Phi^{-1},
\]
where
\[
\sigma_3=
\begin{pmatrix}
1&0\\
0&-1
\end{pmatrix}
\]
is the third Pauli matrix.

Since
\[
SU(2)/U(1)\cong S^2,
\]
the adjusted Gauss map may again be regarded as a sphere-valued map.
Its stereographic projection provides the meromorphic function
appearing in the generalized Weierstrass--Kenmotsu representation.
Although its target coincides with that of the classical Euclidean
Gauss map, its geometric interpretation is fundamentally different.
Rather than representing the unit normal vector itself, it reflects the
balanced normalization of the moving frame that underlies the
representation.

A principal advantage of this construction is that it provides a natural
normalization for generalized Weierstrass representations of constant
mean curvature surfaces with
\[
0\le H<1.
\]
Instead of working directly with the Legendre lift, one considers a
rank-one holomorphic potential
\[
\eta
\in
\Omega^{1,0}(M,\mathfrak{sl}(2,\mathbb C))
\]
together with the associated family of connection forms
\[
\Omega
=
\eta
-
\lambda\eta^{*},
\]
where the spectral parameter satisfies
\[
|\lambda|\le1.
\]
The mean curvature is then recovered through
\[
H
=
\frac{1-|\lambda|^{2}}
{1+|\lambda|^{2}}.
\]
This representation extends the classical philosophy of
Weierstrass-type formulas to the entire range
\(
0\le H<1
\),
providing a unified local construction for these immersions
\cite{TodaGulerAtampalage}.

A fundamental feature of this representation is the appearance of the
connection
\[
\Omega=\eta-\lambda\eta^{*},
\]
whose flatness condition
\[
d\Omega+\Omega\wedge\Omega=0
\]
arises naturally from the balanced matrix formulation and constitutes
the differential integrability condition of the reconstruction process
\cite{TodaGulerAtampalage}. The balancing procedure leads canonically to
this connection, while its vanishing curvature guarantees the
compatibility of the associated linear system and hence the existence of
the corresponding immersion.

It is important to note, however, that this flatness condition is not a
consequence of the algebraic rank-one assumption
\[
\det\eta=0.
\]
The determinant condition guarantees only that the holomorphic
potential has rank one; it does not control the mixed commutator terms
involving \(\eta\) and its Hermitian adjoint. Consequently, the
flatness of \(\Omega\) represents the essential differential
integrability condition underlying the generalized
Weierstrass--Kenmotsu representation.

From a broader perspective, the adjusted Gauss map should therefore be
viewed as a representation-adapted Gauss map rather than merely an
auxiliary construction. Although it is not a canonical invariant of the
immersion, it provides the balanced geometric normalization that links
the moving-frame formulation with explicit reconstruction formulas. It
thus occupies a natural intermediate position between the
contact-geometric Legendre lift and the generalized Weierstrass
representation.

The distinction between these viewpoints is worth emphasizing. The
hyperbolic Gauss maps describe the asymptotic geometry of the normal
congruence, the Legendre lift records the complete contact element of
the immersion, and the adjusted Gauss map provides a balanced
representation of the same geometric data that is particularly well
adapted to explicit reconstruction. These constructions are therefore
complementary rather than competing, each highlighting a different
aspect of hyperbolic surface theory.

The adjusted Gauss map illustrates an important theme recurring
throughout the subject: different geometric questions naturally require
different Gauss-type constructions. The next section presents another
example of this principle by considering the conformal Gauss map, whose
target and geometric interpretation arise from Möbius geometry rather
than from contact geometry or representation theory.
%%%%%%%%%%%%%%%%%%%%%%%%%%%%%%%%%%%%%%%%%%%%%%%%%%%%%%%%%%%%%%%%%%%%%%%%%%%%%%%%%%%%%%%%%%%%%%%%%%%%%%%%%%%%%%%%%%%%%%%%%%%%%%%%%%%%%%%%%%
\section{The Conformal Gauss Map and Willmore Geometry}

The Gauss maps discussed in the previous sections arise naturally from
hyperbolic geometry, contact geometry, and integrable systems. A fourth
and fundamentally different viewpoint comes from conformal differential
geometry.

Rather than associating to an immersed surface the endpoints of its
normal geodesics or a normalized moving frame, conformal geometry
associates to each point its \emph{mean-curvature sphere}. This
construction is invariant under Möbius transformations and therefore
depends only on the conformal structure of the immersion rather than on
a particular realization of the ambient space.

The resulting map is called the \emph{conformal Gauss map}. It assigns
to each point of the surface the unique oriented sphere tangent to the
surface whose mean curvature agrees with that of the immersion. Unlike
the hyperbolic Gauss maps, whose targets lie on the ideal boundary of
hyperbolic space, or the Legendre lift, whose target is the unit tangent
bundle, the conformal Gauss map takes values in the space of oriented
spheres, naturally realized as a homogeneous space of the conformal
group. This viewpoint was developed systematically by Hélein and later
extended by several authors working in integrable systems and Willmore
geometry \cite{Helein,DorfmeisterWang,Wang}.

%%%%%%%%%%%%%%%%%%%%%%%%%%%%%%%%%%%%%%%%%%%%%%%%%%%%%%%%%%%%%%%

\subsection{The Mean-Curvature Sphere Congruence}

For every nonumbilic point of an immersed surface there exists a unique
oriented sphere tangent to the surface whose mean curvature agrees with
that of the immersion. As the point varies over the surface, these
spheres form the \emph{mean-curvature sphere congruence}.

Historically, this family of tangent spheres has been referred to in the
classical differential geometry literature as the
\emph{osculating sphere congruence} or \emph{osculating sphere bundle},
particularly in the French and Eastern European schools. Modern
conformal differential geometry emphasizes its Möbius-invariant
character and therefore generally adopts the terminology
\emph{mean-curvature sphere congruence}, which has become standard in the
theory of conformal Gauss maps and Willmore surfaces.

This sphere congruence provides a conformally invariant analogue of the
classical Gauss map. Whereas the Euclidean Gauss map records only the
oriented normal direction and the hyperbolic Gauss maps record the ideal
endpoints of normal geodesics, the mean-curvature sphere encodes both
first- and second-order geometric information. In particular, it is
preserved under Möbius transformations of the ambient space.

The natural setting for this construction is the light-cone model of
conformal geometry. Identifying the conformal three-sphere with the
projectivized null cone in Lorentzian five-space
\(\mathbb R^{1,4}\), oriented spheres correspond to Lorentzian
four-dimensional subspaces. Consequently, the conformal Gauss map may
be regarded as a map into an appropriate Grassmannian, itself a
Riemannian symmetric space of the conformal group
\cite{Helein}.

Thus, although the construction differs substantially from the Gauss
maps encountered earlier, the underlying philosophy remains the same:
the geometry of an immersed surface is encoded by the analytic
properties of a naturally associated map into a homogeneous target
space.

%%%%%%%%%%%%%%%%%%%%%%%%%%%%%%%%%%%%%%%%%%%%%%%%%%%%%%%%%%%%%%%

\subsection{Hélein's Harmonic Map Formulation}

A major breakthrough was achieved by Hélein, who showed that the
Euler--Lagrange equation for the Willmore functional admits a complete
reformulation in terms of the harmonicity of the conformal Gauss map
\cite{Helein}. More precisely, away from umbilic points, the conformal
Gauss map of an immersion is harmonic if and only if the immersion is a
Willmore surface.

This result places Willmore surface theory squarely within the general
framework of harmonic maps into symmetric spaces. In particular, the
nonlinear fourth-order Willmore equation is replaced by the harmonic-map
equation for the conformal Gauss map, allowing many techniques from
integrable systems to be applied to conformal surface geometry.

The analogy with the classical Ruh--Vilms theorem is striking. In
Euclidean geometry, the harmonicity of the Gauss map characterizes
constant mean curvature surfaces. In conformal geometry, the harmonicity
of the conformal Gauss map characterizes Willmore surfaces. Thus, in
both theories, the Gauss map provides the fundamental geometric object
through which the governing variational equations become harmonic map
equations.

Since harmonic maps into symmetric spaces admit associated families of
flat connections, spectral deformations, and loop-group
factorizations, Hélein's reformulation immediately opens the door to the
full machinery of integrable systems.

%%%%%%%%%%%%%%%%%%%%%%%%%%%%%%%%%%%%%%%%%%%%%%%%%%%%%%%%%%%%%%%

\subsection{Willmore Surfaces and Integrable Systems}

The Willmore functional
\[
\mathcal W(f)
=
\int_M(H^2-K)\,dA
\]
is invariant under conformal transformations of the ambient space.
Critical points of this functional are called \emph{Willmore surfaces}
and occupy a central position in conformal differential geometry.
Earlier variational aspects of this theory are closely related to the
moving-frame methods developed by Bryant and Griffiths
\cite{BryantGriffiths}.

The harmonic-map formulation initiated by Hélein has proved remarkably
fruitful. Xia and Shen constructed generalized Weierstrass-type
representations for Willmore surfaces in higher-dimensional spheres
\cite{XiaShen}. Dorfmeister and Wang subsequently developed a DPW
construction for Willmore surfaces based on the conformal Gauss map
\cite{DorfmeisterWang}, while Wang further clarified the loop-group
approach and its geometric applications \cite{Wang}.

These developments demonstrate that the integrable-systems techniques
originally developed for constant mean curvature surfaces extend
naturally to the broader setting of conformal variational geometry.

%%%%%%%%%%%%%%%%%%%%%%%%%%%%%%%%%%%%%%%%%%%%%%%%%%%%%%%%%%%%%%%

\subsection{Minimal Surfaces and Willmore Geometry}

One of the most remarkable consequences of conformal surface theory is
that minimal surfaces in different space forms become closely related
through conformal geometry.

Classically, stereographic projection identifies minimal surfaces in
\(\mathbb R^3\) with Willmore surfaces in the conformal three-sphere.
Likewise, the conformal compactification
\[
\mathbb H^3
\hookrightarrow
\mathbb S^3
\]
allows every immersed surface in hyperbolic space to be regarded as a
conformal immersion into the three-sphere. Under this compactification,
minimal surfaces in hyperbolic space naturally become Willmore
surfaces in the conformal sense.

From this perspective, Xia and Shen showed that the conformal Gauss map
of a minimal surface in hyperbolic space is harmonic
\cite{XiaShen}. Consequently, minimal surfaces in hyperbolic space admit
two complementary geometric interpretations. On the one hand, they may
be studied through their hyperbolic Gauss maps, Legendre lifts, and
representation formulas. On the other hand, after conformal
compactification they become Willmore surfaces whose geometry is encoded
by the harmonic conformal Gauss map.

This relationship considerably enlarges the unified picture developed
throughout the present survey. Hyperbolic Gauss maps describe the
geometry at infinity; the Legendre lift records the complete contact
element; adjusted Gauss maps provide balanced representations adapted to
explicit Weierstrass constructions; and the conformal Gauss map reveals
the Möbius-invariant geometry visible after conformal compactification.

Rather than representing competing constructions, these Gauss maps
capture complementary aspects of the same immersed surface. Together
they illustrate the rich interplay between hyperbolic geometry,
harmonic maps, integrable systems, and conformal differential geometry.
%%%%%%%%%%%%%%%%%%%%%%%%%%%%%%%%%%%%%%%%%%%%%%%%%%%%%%%%%%%%%%%%%%%%%%%%%%%%%%%

\section{A Unified Perspective on Gauss Maps in Hyperbolic Surface Theory}

The preceding sections have presented four distinct Gauss-type
constructions arising naturally in the theory of immersed surfaces in
hyperbolic three-space. Historically, these maps were introduced at
different times, motivated by different geometric problems, and developed
within different mathematical communities. Consequently, they often
appear in the literature as unrelated constructions.

From the viewpoint developed in this survey, however, these maps should
not be regarded as competing definitions of the same object. Rather, each
one captures a different geometric structure naturally associated with an
immersed surface. Their coexistence reflects not an ambiguity in the
notion of a Gauss map, but the remarkable richness of hyperbolic
geometry itself.

%%%%%%%%%%%%%%%%%%%%%%%%%%%%%%%%%%%%%%%%%%%%%%%%%%%%%%%%%%%%%%%

\subsection{Why Several Gauss Maps?}

In Euclidean three-space, the unit normal vector provides a canonical
choice of Gauss map because translations identify tangent spaces at
different points. Hyperbolic space possesses no comparable global
translation structure. Instead, it carries several independent geometric
structures, each leading naturally to its own Gauss-type construction.

The ideal boundary gives rise to the hyperbolic Gauss maps, which record
the asymptotic behavior of the oriented normal geodesics.

The contact structure on the unit tangent sphere bundle produces the
Legendre Gauss lift, retaining the complete position--normal data of the
immersion.

Gauge normalization and loop-group factorization naturally lead to the
adjusted Gauss maps, which reorganize the same geometric information into
a form adapted to generalized Weierstrass representations.

Finally, conformal compactification introduces the conformal Gauss map,
whose target is the space of mean-curvature spheres and whose harmonicity
characterizes Willmore surfaces.

Each construction therefore answers a different geometric question, and
none should be expected to replace the others.

%%%%%%%%%%%%%%%%%%%%%%%%%%%%%%%%%%%%%%%%%%%%%%%%%%%%%%%%%%%%%%%

\subsection{Comparison of the Principal Gauss Maps}

The principal features of the Gauss maps discussed in this survey are
summarized in Table~\ref{tab:comparison}.

\begin{table}[htbp]
\centering
\caption{Comparison of the principal Gauss-type constructions appearing
in hyperbolic surface theory.}
\label{tab:comparison}

\renewcommand{\arraystretch}{1.25}
\setlength{\tabcolsep}{5pt}
\small

\begin{tabularx}{\textwidth}{
@{}
>{\raggedright\arraybackslash}p{2.6cm}
>{\raggedright\arraybackslash}p{2.8cm}
>{\raggedright\arraybackslash}X
>{\raggedright\arraybackslash}p{3.7cm}
@{}
}
\toprule
\textbf{Gauss map}
&
\textbf{Target space}
&
\textbf{Primary geometric role}
&
\textbf{Distinguished analytic property}
\\
\midrule

Classical Gauss map
&
$\mathbb S^2$
&
Unit normal field
&
Harmonic for CMC surfaces
\\[2pt]

Hyperbolic Gauss maps
&
$\partial_\infty\mathbb H^3$
&
Asymptotic geometry of normal geodesics
&
Conformal for $H=\pm1$
\\[2pt]

Legendre Gauss lift
&
$U\mathbb H^3$
&
Complete contact element
&
Harmonic lift in the generalized DPW framework
\\[2pt]

Adjusted Gauss map
&
Gauge-normalized sphere-valued target
&
Generalized Weierstrass reconstruction
&
Flat rank-one potential
\\[2pt]

Conformal Gauss map
&
Space of oriented mean-curvature spheres
&
Möbius geometry
&
Harmonic for Willmore surfaces
\\

\bottomrule
\end{tabularx}
\end{table}

%%%%%%%%%%%%%%%%%%%%%%%%%%%%%%%%%%%%%%%%%%%%%%%%%%%%%%%%%%%%%%%

\subsection{Relationships Among the Various Constructions}

Although these Gauss maps have different target spaces and different
geometric interpretations, they are closely related. Each is obtained by
emphasizing a particular aspect of the same immersed surface.

The hyperbolic Gauss maps arise by projecting the oriented normal
geodesics to the ideal boundary.

The Legendre lift retains the full position--normal pair and therefore
contains strictly more local information than the hyperbolic Gauss maps.

The adjusted Gauss map is obtained after gauge normalization and Iwasawa
factorization of the moving frame, providing a representation especially
well adapted to explicit reconstruction formulas.

The conformal Gauss map, by contrast, is obtained only after passing to
the conformal compactification of hyperbolic space, where the immersion
is viewed through its mean-curvature sphere congruence rather than its
normal geometry.

These relationships may be summarized schematically as

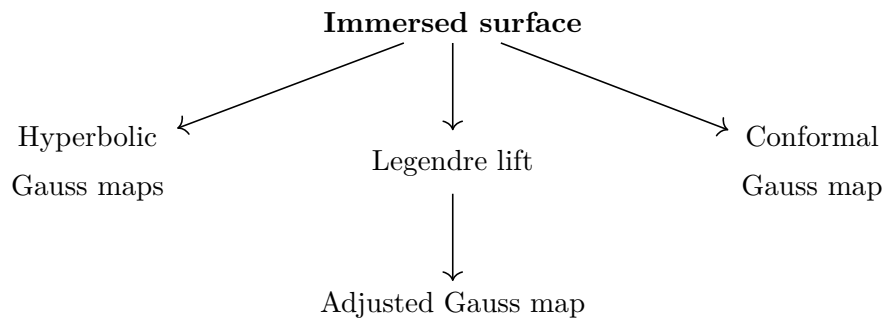
\[
\begin{tikzcd}[
    row sep=2.2em,
    column sep=4.5em,
    arrows={line width=0.6pt}
]
&
\text{\bfseries Immersed surface}
\arrow[dl]
\arrow[d]
\arrow[dr]
&
\\
\begin{gathered}
\text{Hyperbolic}\\
\text{Gauss maps}
\end{gathered}
&
\begin{gathered}
\text{Legendre lift}
\end{gathered}
\arrow[d]
&
\begin{gathered}
\text{Conformal}\\
\text{Gauss map}
\end{gathered}
\\
&
\text{Adjusted Gauss map}
&
\end{tikzcd}
\]

This diagram is not intended to represent canonical morphisms between the
various target spaces. Rather, it illustrates the different geometric
viewpoints through which the same immersed surface may be studied.

%%%%%%%%%%%%%%%%%%%%%%%%%%%%%%%%%%%%%%%%%%%%%%%%%%%%%%%%%%%%%%%

\subsection{Concluding Remarks}

The theory of immersed surfaces in hyperbolic three-space illustrates a
general phenomenon that occurs throughout modern differential geometry:
a single geometric object often admits several complementary
interpretations.

From the asymptotic viewpoint, the hyperbolic Gauss maps describe the
behavior of normal geodesics at infinity. 
From the contact-geometric viewpoint, the Legendre lift provides the
natural framework for harmonic maps into $4$-symmetric spaces and the
generalized DPW method. From the representation-theoretic viewpoint, adjusted Gauss maps lead to explicit generalized Weierstrass--Kenmotsu formulas through gauge
normalization and Iwasawa splitting. Finally, from the conformal viewpoint, the conformal Gauss map reveals
the Möbius-invariant geometry of the immersion and establishes profound
connections with Willmore surfaces and harmonic map theory.

Rather than competing with one another, these constructions complement
each other. Together they provide a unified geometric picture connecting
hyperbolic geometry, harmonic maps, integrable systems, contact
geometry, and conformal differential geometry. We hope that this
perspective clarifies the role of the principal Gauss maps currently
appearing in the literature and provides a useful framework for future
developments in hyperbolic surface theory.

\vspace{2em}
\noindent

\textbf{Magdalena Toda}\\
College of Natural and Applied Sciences, Missouri State University,\\
Springfield, MO 65897, USA\\
\textit{E-mail:}
\texttt{mtoda@missouristate.edu}\\
Department of Mathematics and Statistics, Texas Tech University,\\
Lubbock, TX 79409, USA\\
\textit{E-mail:}
\texttt{magda.toda@ttu.edu}

\vspace{1em}
\textbf{Erhan G\"uler}\\
Department of Mathematics and Statistics, Texas Tech University,\\
Lubbock, TX 79409, USA\\
\textit{Email:} \texttt{eguler@ttu.edu}

\end{document}